\documentclass[11pt]{amsart}

\usepackage{amsmath}
\usepackage{amssymb}

\newtheorem{theorem}{Theorem}[section]
\newtheorem{claim}[theorem]{Claim}

\theoremstyle{definition}
\newtheorem{definition}[theorem]{Definition}

\theoremstyle{remark}

\newcount\skewfactor
\def\mathunderaccent#1#2 {\let\theaccent#1\skewfactor#2
\mathpalette\putaccentunder}
\def\putaccentunder#1#2{\oalign{$#1#2$\crcr\hidewidth
\vbox to.2ex{\hbox{$#1\skew\skewfactor\theaccent{}$}\vss}\hidewidth}}

\def\smallbox#1{\leavevmode\thinspace\hbox{\vrule\vtop{\vbox
   {\hrule\kern1pt\hbox{\vphantom{\tt/}\thinspace{\tt#1}\thinspace}}
   \kern1pt\hrule}\vrule}\thinspace}

\newcommand{\rest}{{\upharpoonright}}

\newcommand{\conc}{{}^\frown\!}

\newcommand{\cf}{{\rm cf}}

\newcommand{\st}{{such that}}

\newcommand{\contr}{{contradiction}}

\newcommand{\then}{{\underline{then}}}
\newcommand{\Then}{{\underline{Then}}}

\newcommand{\mat}{\mathcal}

\title{The Erd\"os-Rado arrow for singular}
\author{Saharon Shelah}
\address{Institute of Mathematics
 The Hebrew University of Jerusalem
 Jerusalem 91904, Israel
 and  Department of Mathematics
 Rutgers University
 New Brunswick, NJ 08854, USA}
\email{shelah@math.huji.ac.il}
\urladdr{http://shelah.logic.at}
\thanks{Research supported by the United States-Israel Binational
Science Foundation. Publication 881}

\begin{document}

\begin{abstract}
We prove that if $\cf(\lambda)>\aleph_0$ 
and $2^{\cf (\lambda)}<\lambda$ then $\lambda \rightarrow 
(\lambda,\omega+1)^2$ in ZFC
\end{abstract}

\maketitle

\section{introduction}

For regular uncountable $\kappa$, the Erd\"os-Dushnik-Miller theorem,
Theorem 11.3 of \cite{EHMR}, states that
$\kappa \rightarrow (\kappa,\ \omega + 1)^2$.
For singular cardinals, $\kappa$, they were only able to obtain
the weaker result, Theorem 11.1 of \cite{EHMR}, that
$\kappa \rightarrow (\kappa,\ \omega)^2$.
It is not hard to see that if $\cf(\kappa) = \omega$ then
$\kappa \not\rightarrow (\kappa,\ \omega + 1)^2$.
If $\cf\ (\kappa) > \omega$ and
$\kappa$ is a strong limit cardinal, then it follows from the
General Canonization Lemma, Lemma 28.1 in \cite{EHMR}, that
$\kappa \rightarrow (\kappa,\ \omega + 1)^2$.  Question 11.4 of \cite{EHMR}
is whether this holds without the assumption that $\kappa$ is a strong limit
cardinal, e.g., whether,
in ZFC,

$$ 
(1)\ \ \aleph_{\omega_1} \rightarrow (\aleph_{\omega_1},\ \omega + 1)^2.
$$

In \cite{ShSt:419} it was proved that 
$\lambda\rightarrow (\lambda,\omega+1)^2$ if 
$2^{\cf(\lambda)}<\lambda$ and there is a nice 
filter on $\kappa$, (see \cite[Ch.V]{Sh:g}: 
follows from suitable failures of SCH). 
Also proved there are consistency results when 
$2^{\cf(\lambda)}\geq \lambda$

Here continuing \cite{ShSt:419} but not relying on it, we 
eliminate the extra assumption, i.e, we prove (in ZFC)

\begin{theorem} 
\label{0.1}
If $\aleph_0 < \kappa = \cf (\lambda)$ and $2^\kappa < \lambda$ 
 \then\ $\lambda \rightarrow
(\lambda,\ \omega + 1)^2$.
\end{theorem}

Before starting the proof, let us recall the well known 
definition: 

\begin{definition}
\label{0.2}
Let $D$ be an $\aleph_1$-complete filter on $Y$, 
and $f\in {}^Y{\rm Ord}$, and $\alpha\in {\rm Ord}
\cup \{\infty\}$.

We define when ${\rm rk}_D (f)=\alpha$
by induction on $\alpha$ (it is well known that 
${\rm rk}_D (f)<\infty$):
\begin{enumerate}
\item[(*)] ${\rm rk}_D (f)=\alpha$
iff $\beta<\alpha\Rightarrow {\rm rk}_D (f)\neq \beta$, 
and for every $g\in {}^Y{\rm Ord}$ satisfying $g<_D f$, 
there is $\beta<\alpha$ \st\ ${\rm rk}_D (g)=\beta$.
\end{enumerate}
Notice that we will use normal filters on $\kappa=
\cf(\kappa)>\aleph_0$, so the demand of $\aleph_1$-
completeness in the definition, holds for us.\newline
Recall also
\end{definition}

\begin{definition}
\label{0.3}
Assume $Y,D,f$ are as in definition \ref{0.2}.
$$
J[f,D]= \{Z\subseteq Y:Y\setminus Z\in D \hbox{ or }
{\rm rk}_{D+(Y\setminus Z)} (f)>{\rm rk}_D (f)\}
$$
Lasly, we quote the next claim
(the definition \ref{0.3} and claim are from \cite{Sh:71},
and explicitly \cite{Sh:589}(5.8(2),5.9)):
\end{definition}

\begin{claim}
\label{0.4}
Assume $\kappa>\aleph_0$ is realized, and
$D$ is a $\kappa$-complete (a normal) filter 
on $Y$.

\Then\ $J[f,D]$ is a $\kappa$-complete (a normal) ideal 
on $Y$ disjoint to $D$ for any $f\in {}^Y{\rm Ord}$
\end{claim}

\section{The proof}

In this section we prove Theorem \ref{0.1} 
of the Introduction, which, for
convenience, we now restate.

\begin{theorem}
\label{1.1}
 If $\aleph_0 < \kappa = \cf(\lambda),\ 2^\kappa < \lambda$ 
\then\ $\lambda \rightarrow
(\lambda,\ \omega + 1)^2$.
\end{theorem}

\begin{proof}
${}$

\par \noindent 
\underline{Stage A}\ 
We know that $\aleph_0 < \kappa = \cf (\lambda) < 
\lambda,\ 2^\kappa < \lambda$ 
We will show that $\lambda \rightarrow (\lambda,\ \omega + 1)^2$.

So, towards a contradiction, suppose that 
\begin{enumerate}
\item[$(*)_1$] $c:[\lambda]^2\ \rightarrow\
\{\hbox{red, green}\}$ but has no red set of cardinality 
$\lambda$ and no
green set of order type $\omega + 1$.  
\end{enumerate}
Choose $\bar \lambda$ \st: 
\begin{enumerate}
\item[$(*)_2$] $\bar \lambda=\langle \lambda_i:i
<\kappa\rangle$ is 
increasing and continuous with limit $\lambda$, 
and for $i = 0$
or $i$ a successor ordinal, $\lambda_i$ is a successor cardinal.
We also let $\Delta_0 = \lambda_0$ and
for $i < \kappa,\ \Delta_{1 + i} = [\lambda_i,\ 
\lambda_{i + 1})$.  For
$\alpha < \lambda$ we will let ${\bf i}(\alpha) = $ the unique 
$i < \kappa$ such
that $\alpha \in \Delta_i$.
\end{enumerate}
We can clearly assume, in
addition, that 
\begin{enumerate}
\item[$(*)_3$] $\lambda_0 > 2^\kappa$, for
$i < \kappa,\ \lambda_{i + 1} \geq \lambda_i^{++}$, and
that each $\Delta_i$ is homogeneously red for $c$.
\end{enumerate}
The last is justified by the Erd\"os-Dushnik-Miller theorem for
$\lambda_{i + 1}$, i.e., as $\lambda_{i+1}\rightarrow 
(\lambda_{i+1},\omega+1)^2$ because 
$\lambda_{i+1}$ is regular.
\medskip

\par \noindent
\underline{Stage B}:\ 
For $0 < i < \kappa$, we define ${\rm Seq}_i$ to be
$\{\langle \alpha_0, ..., 
\alpha_{n-1}\rangle:{\bf i}(\alpha_0) < ... 
< {\bf i}(\alpha_{n-1}) < i\}$.  For $\zeta
\in \Delta_i$ and
$\langle \alpha_0, ..., \alpha_{n-1}\rangle = 
{\bar \alpha} \in {\rm Seq}_i$, 
we say ${\bar \alpha}
\in {\mat T}^\zeta$ iff $\{\alpha_0, ..., \alpha_{n-1},\ \zeta\}$ 
is homogeneously green for $c$.
Note that an infinite $\triangleleft$-increasing branch
in ${\mat T}^\zeta$ 
violates the non-existence of a green set of order type
$\omega + 1$, so, 
\begin{enumerate} 
\item[$(*)_4$] ${\mat T}^\zeta$ is well-founded,
that is we cannot find $\eta_0\triangleleft \eta_1
\triangleleft\ldots \triangleleft \eta_n \triangleleft
\ldots$.
\end{enumerate}

Therefore the following definition of a rank function,
${\rm rk}^\zeta$, on ${\rm Seq}_i$
can be carried out.

If $\eta\in {\rm Seq}_i\setminus T^\zeta$
then ${\rm rk}^\zeta (\eta)=-1$.
We define ${\rm rk}^\zeta:{\rm Seq}_i 
\rightarrow {\rm Ord} \cup \{-1\}$ as follows
by induction on the
ordinal $\xi$, we have
${\rm rk}^\zeta({\bar \alpha}) 
=\xi$ iff for all $\epsilon < \xi, {\rm rk}^\zeta 
({\bar \alpha})$ was not 
defined as $\epsilon$ but there is $\beta$
such that ${\rm rk}^\zeta 
({\bar \alpha}^{\conc}\langle 
\beta\rangle) \geq \epsilon$.  Of course, 
if $\xi$ is a succesor ordinal, it is enough to 
check for $\epsilon=\xi-1$, and for limit
ordinals, $\delta$, if for all $\xi < \delta,\ 
{\rm rk}^\zeta({\bar \alpha}) \geq \xi$,
then ${\rm rk}^\zeta
({\bar \alpha}) \geq \delta$. 
In fact, it is clear
that the range of ${\rm rk}^\zeta$ is a 
proper initial segment of $\mu_i^+$, where
$\mu_i := {\rm card} 
(\bigcup\ \{\Delta_\epsilon:\epsilon < i\})$, 
and so, in particular,
the range of ${\rm rk}^\zeta$ has cardinality 
at most $\lambda_i$.  
Note that $\lambda_{i + 1}
\geq \lambda^{++}_i> \mu_i^+$.

Now we can choose $B_i$, 
an end-segment of $\Delta_i$ such that for all
${\bar \alpha} \in {\rm Seq}_i$ 
and all $0 \leq \gamma < \mu_i^+$,
if there is $\zeta \in B_i$ such that 
${\rm rk}^\zeta({\bar \alpha}) = \gamma$, then
there are $\lambda_{i + 1}$ such $\zeta$-s.  
Recall that $\Delta_{i}$ and therefore
also $B_i$ are of order type $\lambda_{i + 1}$, which is a successor
cardinal $>\mu^+_i> |{\rm Seq}_i|$ hence such $B_i$ exists.  
Everything is now in place for the main definition.
\medskip

\par \noindent 
\underline{Stage C}:\ 
$({\bar \alpha},\ Z,\ D,\ f) \in K$ iff
\begin{enumerate}
\item $D$ is a normal filter on $\kappa$,
\item $f:\kappa \rightarrow {\rm Ord}$,
\item $Z \in D$ 
\item for some $0 < i < \kappa$ we have 
${\bar \alpha} \in {\rm Seq}_i$ and $Z$ is disjoint to $i+1$
and for every
$j \in Z$ (hence $j>i$) there is $\zeta \in B_j$ such that
${\rm rk}^\zeta({\bar \alpha}) 
= f(j)$ (so, in particular, 
${\bar \alpha} \in {\mat T}^\zeta$).
\end{enumerate}
\medskip

\par \noindent
\underline{Stage D}:\ 
Note that $K \neq \emptyset$, since if we choose $\zeta_j \in
B_j$, for $j < \kappa$, take $Z = \kappa\setminus \{0\},\ 
{\bar \alpha} = $
the empty sequence, choose $D$ to be any normal filter on
$\kappa$ and define $f$ by 
$f(j) = {\rm rk}^{\zeta_j}({\bar \alpha})$,
then $({\bar \alpha},\ Z,\ D,\ f) \in K$.

Now clearly by \ref{0.2},
among the quadraples $({\bar \alpha},\ Z,\ D,\ f) \in K$,
there is one with ${\rm rk}_D(f)$ minimal.
So, fix one such quadraple, and denote it by 
$({\bar \alpha}^*,\ Z^*,\ D^*,\ f^*)$.
Let $D^*_1$ be the filter on $\kappa$ 
dual to $J [f^*,D^*]$, 
so by claim \ref{0.4} it is 
a normal filter on $\kappa$ extending $D^*$.

For $j \in Z^*$, set $C_j = \{\zeta \in 
B_j:{\rm rk}^\zeta({\bar \alpha}^*) =
f^*(j)\}$.  Thus 
by the choice of $B_j$ we know that 
${\rm card}(C_j) = \lambda_{j + 1}$, and
for every $\zeta \in C_j$ the set
$({\rm Rang}(\bar \alpha^*) \cup \{\zeta\})$ 
is homogeneously green under the colouring $c$.
Now: suppose $j \in Z^*$.  
For every $\Upsilon \in Z^* \setminus (j + 1)$
and $\zeta \in C_j$,
let $C^+_\Upsilon (\zeta) = \{\xi \in C_\Upsilon:
c(\{\zeta,\xi\}) = 
\hbox{green}\}$.
Also, let $Z^+(\zeta) = \{\Upsilon \in Z^* 
\setminus (j + 1): {\rm card}  
(C^+_\Upsilon (\zeta))
= \lambda_{\Upsilon + 1}\}$.  
\medskip

\par \noindent
\underline{Stage E}:\ 
For $j \in Z^*$ and $\zeta \in C_j$, let
$Y(\zeta) = Z^* \setminus 
Z^+(\zeta)$.  Since $\lambda_0 > 
2^\kappa$ and $\lambda_{j+1}>\lambda_0$ is regular,
for each $j \in Z^*$ there are
$Y = Y_j \subseteq \kappa$ and 
$C'_j \subseteq C_j$ with
${\rm card} (C'_j) = \lambda_{j + 1}$ 
such that $\zeta \in C'_j \Rightarrow
Y(\zeta) = Y_j$.

Let $\hat Z = \{j \in Z^*:Y_j \in D^*_1\}$. 
Now the proof split to two cases. 
\medskip

\par \noindent
\underline{Case 1}:\  
$\hat Z\neq \emptyset\ {\rm mod}\ D^*_1$ 

Define $Y^*=\{j\in \hat Z$: for every $i\in \hat Z 
\cap j$, we have $j\in Y_i\}$. Notice that $Y^*$ is the 
intersection of $\hat Z$ with the diagonal
intersection of $\kappa$ sets from $D^*_1$ 
(since $i\in \hat Z\Rightarrow Y_i\in D^*_1$), 
hence (by the normality of $D^*_1$) $Y^*\neq \emptyset\ 
{\rm mod}\ D^*_1$.
But then, by shrinking the $C'_j$ 
for $j\in Y^*$, 
we can get 
a homogeneous red set of cardinality $\lambda$, 
which is contrary to the assumption toward contradiction.

We define $\hat C_j$ for $j \in Y^*$ by induction on $j$ such
that
$\hat C_j$ is a subset of $C'_j$ of cardinality $\lambda_{j + 1}$.
Now, for $j \in Y^*$, let $\hat C_j$ be the set of
$\xi \in C'_j$
\st\ for every $i\in Y^* \cap j$ and every $\zeta \in \hat C_i$
we have $\xi
\not\in C^+_j (\zeta)$.  So, in fact, $\hat C_j$ 
has cardinality $\lambda_{j+1}$ as it is the result of
removing $<\lambda_{j+1}$ elements from $C'_j$
where $|C'_j|=\lambda_{j+1}$ by its choice. 
That is, the number of such pairs 
$(i,\zeta)$ is $\leq\lambda_j$ and:
for $i\in Y^* \cap j$ and $\zeta\in \hat C_i$:
\begin{enumerate}
\item[(a)] $j\in Y_i$ 
[Why? by the definition of $Y^*$ as $j\in Y^*$]
\item[(b)] $\zeta\in C'_i$ 
[Why? as $\zeta\in \hat C_i$ and $\hat C_i\subseteq 
C'_i$ by the induction hypothesis]
\item[(c)] $Y(\zeta)=Y_i$ 
[Why? as by (b) we have $\zeta\in C'_i$ and the choice 
of $C'_i$]
\item[(d)] $j\in Y(\zeta)$
[Why? by (a)+(c)]
\item[(e)] $j\notin Z^+ (\zeta)$ 
[Why? by (d) and the choice of $Y(\zeta)$
as $Z^*\setminus Z^+ (\zeta)$]
\item[(f)] $C^+_j (\zeta)$ has cardinality 
$<\lambda_{j+1}$
[Why? by (e) and the choice of $Z^+(\zeta)$, as 
$j\in \hat Z \subseteq Z^*$]  
\end{enumerate}
So $\hat C_j$ is a well defined subset of 
$C'_j$ of cardinality $\lambda_{j+1}$ for every 
$j\in Y^*$.
But then, clearly
the union of the $\hat C_j$ for $j \in Y^*$, 
call it $\hat C$ satisfies: 
\begin{enumerate}
\item[$(\alpha)$] it has cardinality $\lambda$ 
[as $j\in Y^*\Rightarrow |\hat C_j|=\lambda_{j+1}$
and ${\rm sup} (Y^*)=\kappa$ as $Y^*\neq \emptyset\  
{\rm mod}\  D^*_1$]
\item[$(\beta)$] $c\rest [\hat C_j]^2$ is
constantly red [as we are assumming $(*)_3$]
\item[$(\gamma)$] if $i<j$ are from $Y^*$ and 
$\zeta\in \hat C_i,\xi\in \hat C_j$ then 
$c\{\zeta,\xi\}={\rm red}$ 
[as $\xi\notin\ C^+_j (\zeta)$]
\end{enumerate}
So $\hat C$ has cardinality $\lambda$ and
is homogeneously red.
This concludes the proof in the case 
$\hat Z \neq \emptyset\ {\rm mod}\ D^*_1$
\medskip

\par \noindent
\underline{Case 2}:
$\hat Z=\emptyset\ {\rm mod}\  D^*_1$. 
In that case there are $i\in Z^*, \beta \in C_i$ 
\st\ $Z^+ (\beta)\neq \emptyset\  
{\rm mod}\  D^*_1$

[Why? well, $Z^*\in D^*\subseteq D^*_1$ and $\hat Z=\emptyset\  
{\rm mod}\ D^*_1$, hence $Z^*\setminus \hat Z\neq \emptyset$.
Choose $i\in Z^* \setminus \hat Z$. By the definition of 
$\hat Z$, $Y_i\notin D^*_1$. So, if $\beta\in C'_i$ 
then $Y(\beta)=Y_i\notin D^*_1$ and choose $\beta\in C'_i$,
so $Y(\beta)\notin D^*_1$ hence by the definition of 
$Y(\beta)$ we have $Z^*\setminus 
Z^+(\beta)=Y(\beta)\notin D^*_1$. Since $Z^*\in D^*_1$, we 
conclude that $Z^+ (\beta) \neq \emptyset\ {\rm mod}\ 
D^*_1$].
Let $\bar \alpha'=\bar \alpha^*\conc\langle \beta\rangle, 
Z'=Z^+ (\beta), D'=D^*+Z'$, it is a normal filter 
by the previous sentence as $D^*\subseteq D^*_1$ and 
lastly we define $f'\in {}^\kappa{\rm Ord}$ 
by: 
\begin{enumerate}
\item[(a)] if $j\in Z'$ then $f' (j)={\rm Min}
\{{\rm rk}^\gamma (\bar \alpha'):
\gamma\in C^+_j(\beta)\subseteq B_j\}$
\item[(b)] otherwise $f'(j)=0$ 
\end{enumerate}
Clearly 
\begin{enumerate}
\item[$(\alpha)$] $(\bar \alpha',Z',D',f')\in K$, 
and 
\item[$(\beta)$] $f'<_{D'} f^*$ 

[Why? as $Z'\in D'$ and if $j\in Z'$ then for some
$\gamma\in C^+_j (\beta)$ we have $f'(j)={\rm rk}^\gamma
(\bar \alpha')={\rm rk}^\gamma 
(\bar \alpha^* \conc\langle\beta\rangle)$ 
which by the definition of ${\rm rk}^\gamma$ is 
$>{\rm rk}^\gamma (\bar \alpha^*)=f^* (j)$,
recalling (a) from stage C.]

hence 
\item[$(\gamma)$] ${\rm rk}_{D'} (f')< {\rm rk}_{D'} (f^*)$

[Why? see Definition \ref{0.2}].
\end{enumerate}
But ${\rm rk}_{D'} (f^*)={\rm rk}_{D^*} (f^*)$ 
as $Z'=Z^+ (\beta)\neq 
\emptyset\  {\rm mod}\ D^*_1$ by the definition of $D^*_1$ 
as extending the filter dual  
to $J[f^*,D^*]$, see Definition \ref{0.3}. 
Hence ${\rm rk}_{D'} (f')
< {\rm rk}_{D^*} (f^*)$, so we get a \contr\ 
to the choice of $(\bar \alpha^*,Z^*,D^*,f^*)$.\newline
Clearly at least one of the two cases holds, 
so we are done.
\end{proof}

\end{document}